\documentclass[3p]{elsarticle}

\usepackage{amssymb}
\usepackage{amsthm}
\usepackage[colorlinks=true,citecolor=blue,linkcolor=blue]{hyperref}
\usepackage{listings}
\journal{Discrete Mathematics}

\newtheorem{prop}{Proposition}

\newtheorem{lem}{Lemma}
\newproof{pf}{Proof}
\newdefinition{defn}{Definition}
\newtheorem*{claim}{Claim}

\newcommand\B{\mathbf{B}}

\renewcommand\vec[1]{\overline{#1}}

\newcommand\Pow{\mathcal{P}}

\newcommand\permApp{\cdot}

\newcommand\Perm{S}

\newcommand\eqdef{\stackrel{\smash{\mathsf{def}}}{=}}

\newcommand\Sec[1]{\mathcal{S}\left(#1\mskip.5mu\right)}

\newcommand\sub\subseteq
\renewcommand\sup\supseteq
\newcommand\sq\sqsubseteq
\renewcommand\equiv\approx
\newcommand\implies\Rightarrow

\newcommand\sttt[1]{{\scriptsize\texttt{#1}}}
\newcommand\be{\[\begin{array}[t]{lllllll}}      \newcommand\ee{\end{array}\]}

\newcommand\nodeUp[1]{
\setbox7=\hbox{\strut$#1$}\wd7=0pt
{\box7 \atop \bullet}}
\newcommand\nodeDown[1]{
\setbox7=\hbox{\strut$#1$}\wd7=0pt
{\bullet \atop \box7}}

\newcommand\vin{\mskip2.7mu\vert\raise1pt\hbox{$\scriptscriptstyle\in$}}

\begin{document}

\begin{frontmatter}  

\title{Some Properties of Inclusions of Multisets\\
       and Contractive Boolean Operators}

\author{Pierre Hyvernat\footnote{This work was partially funded by the French ANR project r\'ecr\'e ANR-11-BS02-0010.}}

\address{Universit\'e de Savoie,            \\
         Laboratoire de Math\'ematiques,    \\
         73376 Le Bourget-du-Lac Cedex,     \\
         France}

\ead{Pierre.Hyvernat@univ-savoie.fr}
\ead[url]{http://lama.univ-savoie.fr/\~{}hyvernat/}

\begin{abstract}
  Consider the following curious puzzle: call an~$n$-tuple~$\vec
  X=(X_1,\dots,X_n)$ of sets smaller than another~$n$-tuple~$\vec Y$ if it has
  fewer \emph{unordered sections}. We show that equivalence classes for this
  preorder are very easy to describe and characterize the preorder in terms of
  the simpler pointwise inclusion and the existence of a special increasing
  boolean operator~$f:\B^n\to\B^n$. We also show that contrary to increasing
  boolean operators, the relevant operators are not finitely generated, which
  might explain why this preorder is not easy to describe concretely.
\end{abstract}

\end{frontmatter}

\section*{Introduction: a puzzle}  
Let~$N$ be a (fixed) set and~$n$ be a (fixed) natural number. We can consider the
following partial order on~$\Pow_{\!\!*}(N)^n$, the collection of~$n$-tuples
of \emph{nonempty} subsets of~$N$:
\[
  \vec X \sub \vec Y
  \quad\eqdef\quad
  \prod_{1\le i\le n} X_i \sub \prod_{1\le i\le n} Y_i
\]
where~$\prod_i X_i$ is the usual cartesian product.
Because we restrict to {nonempty} subsets, this preorder coincide with
pointwise inclusion:
\[
  (X_1,\dots,X_n) \sub (Y_1,\dots,Y_n)
  \quad\iff\quad
  \forall 1\le i\le n,\ X_i \sub Y_i
  \ .
\]
We now consider a commutative version of the cartesian product where instead of
the usual ordered~$n$-tuples, we take ``unordered~$n$-tuples''.
\begin{defn}\label{def:section}
  If~$\vec X = (X_1$, \dots,~$X_n)$ is an~$n$-tuple of nonempty subsets of~$N$,
  define~$\Sec{\vec X}$, the set of \emph{unordered sections} of~$\vec X$, as
  \begin{equation}
    \Sec{\vec X}
    \quad\eqdef\quad
    {\Bigg(\prod_{1\le i\le n} X_i\Bigg)}\Big/{\Perm_n}
    \ ,
  \end{equation}
  where~${\_\,}/{\Perm_n}$ denotes quotienting by the action of the
  symmetric group~$\Perm_n$.
\end{defn}
Strictly speaking,~$\Perm_n$ does not really act on~$\prod_i X_i$ but
on~$N^n$. The notion is well-defined because the orbit of an element
of~$\prod_i X_i$ exists even if not all its elements are themselves
in~$\prod_i X_i$.
From now on, we will drop the adjective ``unordered'' and refer to an element
of~$\Sec{\vec X}$ simply as a section of~$\vec X$.
We now define the preorder~$\sq$ on~$\Pow_{\!\!*}(N)^n$:
\begin{defn}\label{def:order}
  If~$\vec X$ and $\vec Y$ are~$n$-tuples of nonempty subsets of~$N$, we
  define~$\vec X \sq \vec Y$ to mean~$\Sec{\vec X} \sub \Sec{\vec Y}$.
  We write~$\vec X\equiv\vec Y$ for~``$\vec X\sq\vec Y$ and~$\vec Y\sq\vec
  X$'', that is, for~$\Sec{\vec X} = \Sec{\vec Y}$.
\end{defn}
The relation~$\sq$ is only a \emph{preorder} because it is not antisymmetric:~$\big(X_{\sigma(1)}, \dots
, X_{\sigma(n)}\big) \equiv (X_1, \dots , X_n)$ for any permutation~$\sigma$.

\bigbreak
\noindent
The aim of this note is to answer the following questions:
\begin{enumerate}
  \item When do we have~$\vec X \equiv \vec Y$?
  \item What is the relation between $\vec X\sq\vec Y$ and $\vec X\sub\vec Y$?
\end{enumerate}
The problem is subtler than it appears and the first question makes for an
interesting puzzle: while elementary, the proof is more complex than what most
people initially think. Readers are thus encouraged to spend a couple of
minutes playing with the problem before reading on.

\paragraph*{Related notions}
The notion of \emph{system of representatives} was introduced by P.~Hall in
1935~\cite{Hall}. A system of representatives for the~$n$-tuple of sets~$\vec
X$ is simply an~$n$-tuple~$\vec x$ such that there is a permutation~$\sigma$
satisfying~$x_i\in X_{\sigma(i)}$ for each~$i\in\{1,\dots,n\}$. Equivalence
classes of those under permutations are exactly the unordered sections
of~$\vec X$ of Definition~\ref{def:section}.
A lot of attention has been devoted to systems of \emph{distinct}
representatives, also called \emph{transversals}, where the components of~$\vec
x$ are pairwise distinct~\cite{Mirsky}.
%
Rather than looking at them individually, we look here at the collection of
all possible systems of representatives. This shift of focus seems to be new
in itself, as is the notion of contractive increasing boolean operator that
appears later. Relating the two will give a concise answer to the second
question.

\smallbreak
This work can also be seen as a first step toward a \emph{factor theory} for
commutative regular algebra~\cite{Conway}. In his book on regular languages,
John Conway develops a fascinating theory of factorization:
if~$R$ is a regular set, a subfactorization is tuple of sets~$\vec
X$ of words satisfying $X_1 \cdot \dots \cdot X_n \sub R$,
where~$\_\cdot\_$ denotes concatenation of regular sets. A
subfactorization is a \emph{factorization} if each~$X_i$ is maximal and a
\emph{factor} of~$R$ is any such~$X_i$. Conway shows in particular that a
regular set has only finitely many factors, and that they are all regular.

Conway devotes a chapter to commutative regular algebra, i.e. the theory
arising from regular algebra when word concatenation is made commutative.
Factor theory isn't part of this chapter, probably because
\textsl{``Commutative regular algebra is notable for the number of results
whose proofs one would expect to be trivial, but which turn out to be very
subtle.''} (\cite{Conway}, page~95). Commutative factor theory certainly looks
very subtle and this work only gives a very partial answer: given the regular
set~$Y_1 \cdot \dots \cdot Y_n$ where each~$Y_i$ is a set of \emph{symbols},
we characterize it factorizations consisting of exactly~$n$ factors.

\smallbreak
The initial motivation for this work comes from a very different area:
denotational models of linear logic. In~\cite{PT}, the relation~$\Sec{\vec
X}\sub T$ played an important role, where the set~$T$ was an arbitrary
collection of~$n$-multisets. Understanding this relation was necessary to
compute small examples, and the preorder~``$\sq$'' naturally appeared in this
way. (Note that the results of this paper are not to make those computation
any easier than they were...)

\paragraph*{Notation} To make formulas less verbose, we will abuse the vector
notation by lifting~``$\in$'' pointwise: just as~$\vec X\sub\vec Y$ means
``$\forall 1\le i\le n,\ X_i\sub Y_i$'', the notation~$\vec a\in\vec X$ is a
synonym for~``$a_i\in X_i$ for all~$1\le i\le n$''.
The (left) action of~$\Perm_n$ on~$n$-tuples is written with a dot and is defined as
$\sigma\permApp\vec a \eqdef\big(a_{\sigma^{-1}(1)}, \dots, a_{\sigma^{-1}(n)}\big)$.
When talking about~$n$-multisets (orbits for the action of~$\Perm_n$), we
identify an~$n$-tuple with its orbit. In particular,~$\vec a \in \Sec{\vec X}$
means that~$\sigma\permApp\vec a\in\vec X$ for some permutation~$\sigma$,
i.e., that~$a_i \in X_{\sigma(i)}$ for all~$1\leq i\leq n$.


\section{The equivalence relation} 
The first question has a simple answer: equivalence is just equality up to a
permutation of the sets. In other words, the failure of antisymmetry is
captured by the remark coming after Definition~\ref{def:order}.
\begin{prop}\label{prop:equivalence}
  Given any~$\vec X$ and~$\vec Y$ in~$\Pow_{\!\!*}(N)^n$, we have
  \begin{equation}\label{eqn:equivalence}
    \vec X\equiv\vec Y
    \quad \iff \quad
    \exists \sigma\in\Perm_n,\ %
    \sigma\permApp\vec X = \vec Y
    \ .
  \end{equation}
\end{prop}
This proposition is slightly surprising because the left side is definitionally equal
to
\begin{equation}\label{eqn:complexe_equivalence}
  \forall \vec a \in \vec X,\ %
  \exists \sigma\in\Perm_n,\ %
  \sigma\permApp\vec a \in \vec Y%
  \quad\hbox{and}\quad
  \forall \vec b \in \vec Y,\ %
  \exists \sigma'\in\Perm_n,\ %
  \sigma'\permApp\vec b \in \vec X%
  \ ,
\end{equation}
while the right side is definitionally equal to
\[
  \exists \sigma\in\Perm_n,\quad%
  \Big(\forall \vec a \in \vec X,\ %
  \sigma\permApp\vec a \in \vec Y%
  \ \hbox{and}\ %
  \forall \vec b \in \vec Y,\ %
  \sigma^{-1}\permApp\vec b \in \vec X\Big)%
  \ .
\]
That the latter implies the former is trivial.
Proposition~\ref{prop:equivalence} asserts the converse:
in~(\ref{eqn:complexe_equivalence}), we can choose the permutation uniformly
for all the~$\vec a\in\vec X$ and~$\vec b\in\vec Y$!

\begin{lem}\label{lem:first}
  We have
  \begin{enumerate}
    \item If~$\vec X\sq\vec Y$ then, for all\/~$1\le j\le n$, there is
      some~$1\le i\le n$ s.t. $X_i \sub Y_{j}$.
    \item If~$\vec X\equiv \vec Y$ then~$X_{i_0}=Y_{j_0}$ for some pair~$i_0,j_0$.
  \end{enumerate}
\end{lem}

\begin{pf}
For the first point, suppose that there is some~$j_0$ satisfying~$X_i \not\sub
Y_{j_0}$ for all~$1\leq i\leq n$. This means that there is an~$\vec a\in\vec
X$ s.t.~$a_i\notin Y_{j_0}$ for all~$i$. This~$\vec a$ cannot be a section
of~$\vec Y$. Contradiction!

The second point follows easily: starting from~$Y_1$ and repeatedly using the
first point, we can construct an infinite sequence~$i_1, j_2, i_3, j_4, \dots$
satisfying:
\[
  \cdots \quad \sub X_{i_{2k+1}} \sub Y_{i_{2k}} \sub \cdots \sub X_{i_{3}}\sub Y_{j_2} \sub X_{i_1} \sub Y_1
\]
Because there are only finitely many possible indices, there are~$k$ and~$k'$,
with~$k<k'$ and~$i_{2k} = i_{2k'}$.
This implies that the sets~$X_{i_{2k}}$ and~$Y_{j_{2k+1}}$ are equal.
\qed
\end{pf}
Thus, if~$\vec X\equiv\vec Y$, one of the sets appears on both sides and we
can start the construction of~$\sigma$ in~$(\ref{eqn:equivalence})$. To finish
the proof of Proposition~\ref{prop:equivalence} by induction on~$n$, we need
to show the following implication:
\begin{equation}\label{eqn:division1}
  (Z,X_2,\dots,X_n) \equiv (Z,Y_2,\dots,Y_n)
  \ \implies\ %
  (X_2,\dots,X_n) \equiv (Y_2,\dots,Y_n)
  \ .
\end{equation}

If the collection of unordered sections is seen as a ``commutative cartesian
product'', the next definition would be the corresponding ``division''.
\begin{defn}
  Let $T$ be a collection of~$n$-multisets and~$Z\in\Pow_{\!\!*}(N)$ we put
\[
  T\div Z
  \quad \eqdef \quad
  \Big\{ (a_2,\dots,a_n)   \ |\ \forall a\in Z,\ (a,a_2,\dots,a_n) \in T \Big\}
  \ .
\]
\end{defn}
When~$Z = \{a\}$, it is a commutative version of Brozozowski's
\emph{derivative}~\cite{broz}, and in the general case, it corresponds to the
commutative notion of \emph{factor} of~$T$ with respect
to~$Z$~\cite{quotients}.

Implication~$(\ref{eqn:division1})$ above follows from the following lemma:
\begin{lem}\label{lem:division}
We have:
\begin{equation}
  \Sec{X_1,X_2,\dots,X_n} \div X_1
  \quad = \quad
  \Sec{X_2,\dots,X_n}
  \ .
\end{equation}
\end{lem}

\begin{pf}
The~``$\sup$'' inclusion follows from the definition.

\smallbreak
For the~``$\sub$'' inclusion, suppose~$(a_2,\dots,a_n) \in \Sec{\vec X} \div
X_1$ and choose~$b\in X_1$. By hypothesis, we necessarily have~$(b,a_2, \dots, a_n) \in
\Sec{\vec X}$, i.e., there is a permutation~$\tau$ s.t.~$(b,a_2, \dots, a_n)
\in (X_{\tau(1)}, \dots, X_{\tau(n)})$.

\smallbreak

If~$\tau(1)=1$, then~$\tau$ defines a permutation on~$\{2, \dots, n\}$ and we
have~$(a_2, \dots, a_n) \in (X_{\tau(2)}, \dots, X_{\tau(n)})$. We can
conclude directly.

\smallbreak

If~$\tau(1)\neq1$, up to permuting the sets~$X_2$, \dots, $X_n$ and choosing
an appropriate element in the orbit of~$(a_2,\dots a_n)$, we can assume
that~$\tau(1)=2$,~$\tau(2)=1$ and~$\tau(i)=i$ when~$2<i\le n$, or in other words,
that~$b\in X_2$, $a_2\in X_1$ and $a_i\in X_i$ whenever~$2<i\le n$.

\noindent
Put~$a_1\eqdef a_2$. Because~$a_1=a_2\in X_1$, we have~$(a_1,a_2, \dots,
a_n)\in\Sec{\vec X}$ by hypothesis, that is,~$\sigma\permApp\vec a \in\vec
X$ for some permutation~$\sigma$. Note that since~$a_1=a_2$, we can
interchange the values~$\sigma(1)$ and~$\sigma(2)$ and still
have~$\sigma\permApp\vec a\in\vec X$.

\noindent
Let~$k \eqdef \min \big\{ i\ |\ i>0, \sigma^i(1)\in \{1,2\}\big\}$. Up to
changing the values of~$\sigma(1)$ and~$\sigma(2)$, we can assume
that~$\sigma^k(1)=1$ and that the set~$I \eqdef
\{1,\sigma(1),\dots,\sigma^{k-1}(1)\}$ is the cycle containing~$1$. We define
the set~$I^{c}$ as~$\{1,\dots,n\}\setminus I$.

\noindent
Rearrange the columns of
\[\vbox{\halign{&$\mskip25mu#\mskip25mu$\hfil\cr
\nodeUp{a_1} & \nodeUp{a_2} & \dots & \nodeUp{a_i} & \dots & \nodeUp{a_n}  \cr
\vin         & \vin         &       & \vin         &       & \vin          \cr
\nodeDown{X_{\sigma(1)}} & \nodeDown{X_{\sigma(2)}} & \dots & \nodeDown{X_{\sigma(i)}} & \dots
& \nodeDown{X_{\sigma(n)}} \cr
}}\]
into two parts:
\[
\underbrace{\vbox{\halign{&$\mskip25mu#\mskip25mu$\hfil\cr
 \nodeUp{a_1}             & \nodeUp{a_{\sigma(1)}}             & \dots\quad & \nodeUp{a_{\sigma^{k-1}(1)}}   \cr
 \vin                     & \vin                     &            & \vin                           \cr
 \nodeDown{X_{\sigma(1)}} & \nodeDown{X_{\sigma^2(1)}} & \dots\quad & \nodeDown{X_{\sigma^k(1)}=X_1} \cr
}}\quad}_{I}
\qquad
\underbrace{\vbox{\halign{&$\mskip25mu#\mskip25mu$\hfil\cr
 \nodeUp{a_2}             & \dots & \nodeUp{a_i}             & \dots \cr
 \vin                     &       & \vin                     &       \cr
 \nodeDown{X_{\sigma(2)}} & \dots & \nodeDown{X_{\sigma(i)}} & \dots\cr
}}}_{{I^{c}}}
\ .
\]
The indices of~$\vec X$ on the left are exactly those in~$I$, and so are the
indices of~$\vec a$. Thus, the indices of~$\vec X$ and~$\vec a$ on the right
are exactly those in~$I^c$.
This shows that~$(a_i)_{i\in{I^{c}}}$ is a section of~$(X_i)_{i\in{I^{c}}}$.
Also, because each of~$\sigma(1)$, \dots, $\sigma^{k-1}(1)$ is strictly
more than~$2$ (by the definition of~$k$), we have~$a_{\sigma^i(1)}\in
X_{\sigma^i(1)}$ for all~$1\le i\le k-1$ by a previous hypothesis.
This shows that the permutation
\[
  \rho : \{2,\dots,n\} \to \{2,\dots,n\},
  \quad \rho(i) \eqdef\cases{
    i & if $i\in\{\sigma(1)$, \dots, $\sigma^{k-1}(1)\}$ \cr
    \sigma(i) & otherwise}
\]
satisfies~$\rho\permApp(a_2, \dots, a_n) \in (X_2, \dots, X_n)$. This
finishes the proof that~$(a_2, \dots, a_n)$ is indeed a section of~$(X_2,
\dots, X_n)$.
\qed
\end{pf}


\section{The preorder}  
The initial question was not very formal and read as: ``\emph{What is the
relation between~$\vec X\sq\vec Y$ and~$\vec X\sub\vec Y$?\/}'' It is obvious
that~$\vec X\sub\vec Y$ implies~$\vec X\sq\vec Y$, but unfortunately, the
converse does not hold, even if we consider~$n$-tuples of sets
up to permutations. For example, we have
\[
  \vec X := \big(\{3\},\{1,2,3\}\big)
  \quad \sq \quad
  \vec Y := \big(\{1,3\},\{2,3\}\big)
\]
because the sections of~$\vec X$ are all sections of~$\vec Y$:
\[
  \Sec{\vec X} = \big\{[3,1],[3,2],[3,3]\big\}
  \quad\subset\quad
  \Sec{\vec Y} = \big\{[1,2],[1,3],[3,2],[3,3]\big\}
  \ .
\]
Lemma~\ref{lem:first} asserts that each set on the right is a superset
of some set on the left. This is indeed the case as both~$\{1,3\}$
and~$\{2,3\}$ are supersets of the same set~$\{3\}$.
However, one set on the left side is strictly bigger than all the sets
on the right side: $\{1,2,3\}\supset\{2,3\}$
and~$\{1,2,3\}\supset\{1,3\}$!

More generally,~$\big(Y_1\cap Y_2 , Y_1\cup Y_2)\sq \big(Y_1,Y_2\big)$ and any
operator~$F$ on~$n$-tuples of sets obtained by composing functions~$(Y_i,Y_j)
\mapsto (Y_i\cap Y_j, Y_i\cup Y_j)$ on any pairs of coordinates,\footnote{provided the
intersection isn't empty to agree with Definition~\ref{def:section}} will
satisfy~$F(\vec Y)\sq \vec Y$. For example,
\be
  F(Y_1,Y_2,Y_3) &\quad\eqdef\quad&
  \Big(Y_1\cap Y_3, Y_2\cap(Y_1\cup Y_3) , Y_2 \cup (Y_1\cup Y_3)\Big)\cr
  &\quad\sq\quad&
  \big(Y_1\cap Y_3,Y_2,Y_1\cup Y_3\big) \cr
  &\quad\sq\quad&
  \big(Y_1,Y_2,Y_3\big)
  \ .
\ee
We will characterize (Proposition~\ref{prop:characterization}) which
operators~$F$ on~$\Pow_{\!\!*}(N)^n$ satisfy~$F(\vec{Y})\sq \vec{Y}$ by
looking at functions acting on tuples of booleans, i.e., \emph{boolean
operators}.

\begin{defn}
  Let $\B \eqdef \{0,1\}$ equipped with the order~$0\le1$. This is a complete
  lattice with operations written~$\vee$ and~$\wedge$. The lattice structure
  is lifted pointwise to~$\B^n$.

\noindent
  If~$u\in\B^n$, the \emph{weight} of~$u$ is the number of~$1$s in~$u$. It is
  written~$|u|$.
\end{defn}
We write elements of~$\B^n$ as words: for example,~``$011101\in\B^6$''.

\begin{defn}
  If~$a\in N$ and~$\vec X\in\Pow_{\!\!*}(N)^n$, the \emph{characteristic
  word of~$a$ along~$\vec{X}$} is an element of~$\B^n$. It is
  written~$\chi_{\vec X}(a)$ and is
  defined by~$\big(\chi_{\vec X}(a)\big)_i = 1$ iff~$a\in X_i$.
\end{defn}
Thus, $\chi_{\vec X}(a)$ describes in which components of~$\vec X$ the
element~$a$ appears, and~$|\chi_{\vec X}(a)|$ is the number of the
components of~$\vec X$ which contain~$a$. There is a necessary condition
for~$\vec X\sq \vec Y$: ``for all~$a\in N$, if~$a$ appears in~$k$ components
of~$\vec Y$, then it appears in at most~$k$ components of~$\vec X$''.
Concisely, this condition can be written as~``$|\chi_{\vec X}(a)| \leq
|\chi_{\vec Y}(a)|$ for all~$a\in N$''. This condition is reminiscent of the
condition appearing in Hall's celebrated ``marriage theorem''~\cite{Hall}.
Like in the marriage theorem, this condition is also sufficient in the
appropriate setting:
\begin{prop}\label{prop:variantHall}
  Given~$\vec X$ and~$\vec Y$ two~$n$-tuples of non-empty subsets of~$N$ that
  satisfy
  \begin{enumerate}
    \item the function~$a\mapsto\chi_{\vec Y}(a)$ is bijective from~$N$
      to~$\B^n\setminus\{0\cdots0\}$,
    \item the function~$f:\chi_{\vec Y}(a)\mapsto\chi_{\vec X}(a)$, with
      domain~$\B^n\setminus\{0\cdots0\}$ is increasing;
  \end{enumerate}
  we have~$\vec X \sq \vec Y$ if and only if $\big|\chi_{\vec X}(a)\big| \le
  \big|\chi_{\vec Y}(a)\big|$ for all~$a\in N$.
\end{prop}
Looking at the example~$\big(\{3\},\{1,2,3\}\big)\sq\big(\{2,3\},\{1,3\}\big)$
might help to understand the conditions of the proposition. The first condition
means that there is exactly one element that belongs only to~$Y_1$ (``$2$''),
exactly one element that belongs only to~$Y_2$ (``$1$'') and exactly one
element that belongs to both (``$3$''). When the first condition is satisfied,
the second condition amounts to ``when~$a$ appears in more sets than~$b$ on
the right side, then~$a$ appears in more~sets than~$b$ on the left side''. The
graph of the resulting function~$f$ can be read below
\be
a : \qquad                 & \quad1\quad   & \quad2\quad   & \quad3\quad   \cr
\chi_{\vec Y}(a) : \qquad  & \quad01\quad  & \quad10\quad  & \quad11\quad  \cr
\chi_{\vec X}(a) : \qquad  & \quad01\quad  & \quad01\quad  & \quad11\quad & \ .\cr
\ee
We can extend this graph with a harmless~$0\cdots0 \mapsto 0\cdots0$ to obtain
the function~$(b_1,b_2) \mapsto (b_1\wedge b_2, b_1\vee b_2)$.
Proposition~\ref{prop:variantHall} follows from a more general lemma:
\begin{lem}\label{lem:orderCharacterization}
We have $\vec X\sq \vec Y$ iff $ f(u) \eqdef \bigvee_{\chi_{\vec Y}(a) \leq u}
\chi_{\vec X}(a) $ satisfies~$|f(u)| \leq |u|$ for all~$u$.
\end{lem}
The function~$f$ is the least increasing function (for the extensional order)
satisfying~$f\big(\chi_{\vec Y}(a)\big)\ge \chi_{\vec X}(a)$ for any~$a$.

\begin{pf}
For the~``$\Leftarrow$'' implication, suppose that~$\vec a\in\vec X$. We want
to show that~$\sigma\permApp\vec a\in\vec Y$ for some permutation~$\sigma$. If
we look at the bipartite graph~$G_{\vec{a},\vec{Y}}$
\[
 G_{\vec{a},\vec{Y}} \quad\eqdef\quad
 \vcenter{\halign{&$\mskip20mu#\mskip20mu$\hfil\cr
  \nodeUp{a_1} & \nodeUp{a_2} & \dots & \nodeUp{a_n} \cr
 \cr
  \nodeDown{Y_1} & \nodeDown{Y_2}  & \dots & \nodeDown{Y_n} \cr
}}\ ,\]
with an edge between~$a_i$ and~$Y_j$ when~$a_i\in Y_j$, finding a~$\sigma$
s.t.~$\sigma\permApp\vec a\in\vec Y$ is equivalent to finding a perfect
matching in~$G_{\vec{a},\vec{Y}}$. By Hall's marriage theorem, this
is equivalent to ``every subset of~$\{a_1, \dots, a_n\}$ of cardinality~$p$ has at
least~$p$ neighbors''.

Take some subset~$U\sub\{a_1, \dots, a_n\}$ of cardinality~$p$.
Because~$\vec a$ is a section of~$\vec X$, this set has at least~$p$ neighbors
in the corresponding~$G_{\vec a,\vec X}$ graph.
Let~$u\eqdef\bigvee_{a\in U} \chi_{\vec Y}(a)$, if~$a\in U$, then~$\chi_{\vec Y}
(a)\leq u$ by definition of~$u$, so that~$\chi_{\vec X}(a) \leq f(u)$ by
definition of~$f$. We thus have~$\bigvee_{a\in U} \chi_{\vec X}(a)\leq f(u)$.
We get
\[
  p
   \quad\le\quad
  \underbrace{\left|\bigvee_{a\in U} \chi_{\vec X}
  (a)\right|}_{\hbox{\footnotesize\# of
  neighboors}\atop\hbox{\footnotesize of $U$ in $G_{\vec a,\vec X}$.}}
  \quad\leq\quad
  | f(u) |
  \quad\leq\quad
  |u|
  \quad\eqdef\quad
  \underbrace{\left| \bigvee_{a\in U} \chi_{\vec Y}(a)\right|}_{\hbox{\footnotesize\# of
  neighboors}\atop\hbox{\footnotesize of $U$ in $G_{\vec a,\vec Y}$.}}
  \ ,
\]
which concludes the~``$\Leftarrow$'' implication.

\smallbreak
For the~``$\implies$'' implication, let~$\vec X\sq\vec Y$ and~$p=|u|<|f(u)|$.
By definition of~$f$, we can find a set~$\{a_1, \dots, a_k\}\sub N$ which
satisfies~$\chi_{\vec Y}(a_i)\le u$ for~$i=1,\dots,k$ and~$\big|\bigvee_{i\le
k} \chi_{\vec X}(a_i)\big|>p$. In particular, we have
\[
  \left|\bigvee_{1\le i\le k} \chi_{\vec Y}(a_i)\right|
  \quad\leq\quad
  p
  \quad<\quad
  \left|\bigvee_{1\le i\le k}\chi_{\vec X}(a_i)\right|
  \ .
\]
For each~$1$ in~$\bigvee_{i\le k}\chi_{\vec X}(a_i)$ take one element
of~$\{a_1, \dots, a_k\}$ that accounts for this~$1$. Call the resulting
tuple~$\vec a$.
Note that this might not be an~$n$-tuple, but its length is strictly greater
than~$p$ (and may contain repetitions). It is only a partial section of~$\vec
X$: to complete it into a section of the whole~$\vec X$, simply add one
element from each of the remaining (non-empty) sets.
The result is also a section of~$\vec Y$ and in particular, each element
of~$\vec a$ needs to fit in one component of~$\vec Y$. This is impossible
because there are at most~$p$ sets~$Y_j$ that can contain the elements of the
tuple~$\vec a$. Contradiction!

\qed
\end{pf}
Lemma~\ref{lem:orderCharacterization} does characterize the~$\sq$ preorder but
still looks a little ad-hoc. We now give a more concise characterization that
relates~$\sq$ with~$\sub$, thus answering our initial question. First
note that we can lift any~$f:\B^n\to\B^m$ to a
function~$\Pow(N)^n\to\Pow(N)^m$:
\begin{defn}
Suppose~$f:\B^n\to\B^m$, define~$\widehat f : \Pow(N)^n\to\Pow(N)^m$ as
\[
  \widehat f\left(\vec Y\right)
  \eqdef
  \vec X
  \quad
  \hbox{with $a\in X_i$ iff $f\big(\chi_{\vec Y}(a)\big)$ has a~$1$ at
  coordinate~$i$}
  \ .
\]
\end{defn}
This transformation is, in a precise categorical sense, \emph{natural}. It
amounts to lifting the boolean operations~$\wedge$ and~$\vee$ to their set
theoretic versions~$\cap$ and~$\cup$ in a way that is compatible with function
composition. For example, with the ``and/or''
function~$(b_1,b_2)\mapsto(b_1\wedge b_2,b_1\vee b_2)$ we obtain~$(Y_1,Y_2)
\mapsto (Y_1\cap Y_2, Y_1\cup Y_2)$.

\begin{defn}
  Call an increasing function~$f:\B^n \to \B^n$ \emph{contractive} if it
  satisfies~$|f(u)| \leq |u|$ for all~$u\in\B^n$.
\end{defn}
A corollary to Lemma~\ref{lem:orderCharacterization} is:

\begin{prop}\label{prop:characterization}
For any~$\vec X$ and~$\vec Y$, we have
\[
  \vec X \sq \vec Y
  \quad\iff\quad
  \vec X \sub \widehat f\left(\vec Y\right)
  \quad\hbox{for some increasing, contractive~$f:\B^n\to\B^n$.}
\]
\end{prop}

\begin{pf}
We know that~$\vec X\sq\vec Y$ is equivalent to having~$|f(u)|\le|u|$ for
all~$u$ in~$\B^n$, where~$f$ is defined as in
Lemma~\ref{lem:orderCharacterization}. This function~$f$ satisfies~$\vec X
\sub \widehat f(\vec Y)$: use~$u\eqdef\chi_{\vec Y}(a)$ to check that~$a\in
X_i$ is an element of the~$i$-th set of~$\widehat f(\vec Y)$.

For the converse, suppose~$f$ is contractive increasing with~$\vec X\sub\widehat
f(\vec Y)$ and let~$\chi_{\vec Y}(a)\le u$. Suppose that~$\chi_{\vec X}(a)$
contains a~$1$ in position~$i$. This means that~$a\in X_i$ and thus~$a$ is
in the~$i$-th set of~$\widehat f(\vec Y)$. We can conclude that~$f\big(\chi_{\vec Y}(a)\big)$
contains a~$1$ in position~$i$. This implies that~$f(u)$ also contains
a~$1$ in position~$i$.

\qed
\end{pf}

\section{Contractive functions are not finitely generated} 
\label{sec:notFinitelyGenerated}
\begin{figure}[ht]
  \centering
  \includegraphics[width=.8\textwidth]{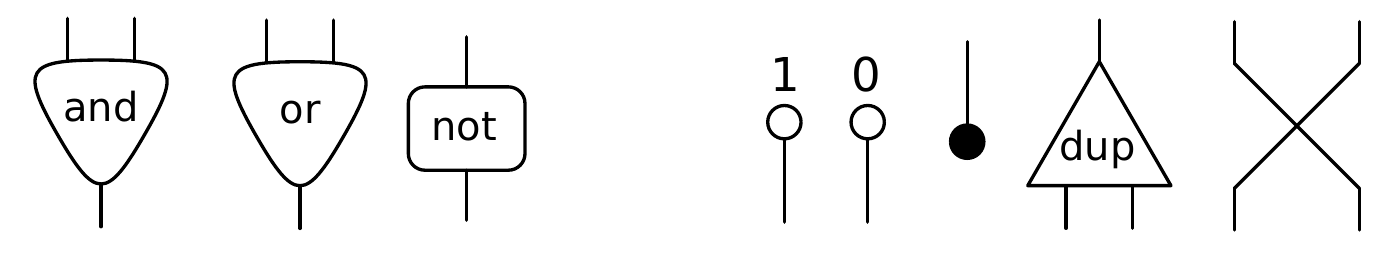}
  \caption{boolean cells}\label{fig:cells}
\end{figure}

If one had a simple representation of contractive increasing boolean operators
from~$\B^n$ to itself, then Proposition~\ref{prop:characterization} would give a simple representation
of the~$\sq$ preorder.
It is well known that all boolean operators~$\B^n\to\B^m$ with~$n$ inputs
and~$m$ outputs can be represented by a boolean circuit using only ``and'',
``or'' together with ``not'' cells. Strictly speaking, we also need constant
values and need a way to forget or duplicate inputs.
The complete set of cells is depicted in Figure~\ref{fig:cells}, where the last cells are:
\begin{itemize}
  \item constants~$1$ and~$0$ (zero input, one output),
  \item drop (one input, zero output),
  \item duplicate (one input, two outputs),
  \item crossing (two inputs, two outputs).
\end{itemize}

These cells, together with a finite set of relations expressing properties of
the operations (associativity, etc.), give a finite presentation of the
{monoidal category} of boolean operators.\footnote{All of this has a precise
algebraic meaning, see~\cite{Lafont} for details.} We can generate the
subcategory of increasing operators by removing the ``not'' cell from the
generators. Unfortunately, no such thing is possible for \emph{contractive
increasing} boolean operators.

\begin{prop}\label{prop:not_finitely_generated}
Contractive increasing boolean operators are not finitely generated.
\end{prop}
In other words, any finite set of cells will either miss some contractive
boolean operator, or generate some non contractive boolean operator.

First, a preliminary lemma:
\begin{lem}\label{lem:bijectionStrictlyIncreasing}
  Let~$f:\B^n \to \B^n$ be an increasing contractive boolean operator; the
  following are equivalent:
  \begin{enumerate}
    \item $f$ is the action of a permutation~$u \mapsto \sigma\permApp u$ for
      some~$\sigma\in\Perm_n$,
    \item $f$ is bijective,
    \item $f$ is injective on words of weight~$1$.
  \end{enumerate}
\end{lem}
\begin{pf}
  Trivially,~$1$ implies~$2$ and~$2$ implies~$3$. Suppose now that~$f$ is
  increasing and contractive on~$\B^n$. Suppose moreover that~$f$ is injective
  on words of weight~$1$. We can define
  a permutation~$\sigma$ on~$\{1, \dots, n\}$ by putting
  \[
  \tau(i) = j
  \quad\hbox{iff}\quad
  f(e_i) = e_j
  \]
  where~$e_i$ represents the word with a single~$1$ in position~$i$. If~$u$
  contains a~$1$ in position~$i$, then, because~$f$ is increasing,~$f(u)$ must
  contain a~$1$ in position~$\tau(i)$. Because~$f$ is contractive,~$f(u)$ cannot contain
  more~$1$s than there are in~$u$. Thus, the~$1$s of~$f(u)$ correspond exactly
  to the images of the~$1$s of~$u$ along~$\tau$: $f$ is indeed the
  action of a permutation.
  \qed
\end{pf}

\begin{pf}[Proposition~\ref{prop:not_finitely_generated}]
Suppose, by contradiction, that there is a finite set of cells that generates
all increasing contractive boolean operator, and write~$m$ for the maximal
arity of the cells in this set.

\noindent
Any non-invertible function has a representation as in
Figure~\ref{fig:normalFormCircuit} where
\begin{itemize}
  \item the topmost rectangle contains only crossings (and invertible cells
    which are, by Lemma~\ref{lem:bijectionStrictlyIncreasing}, equivalent to
    crossings),
  \item the cell~$C$ is not invertible and has arity~$c\leq m$,
  \item and the lowermost rectangle contains the rest of the circuit.
\end{itemize}
\begin{figure}[ht]
  \centering
  \includegraphics[height=5cm]{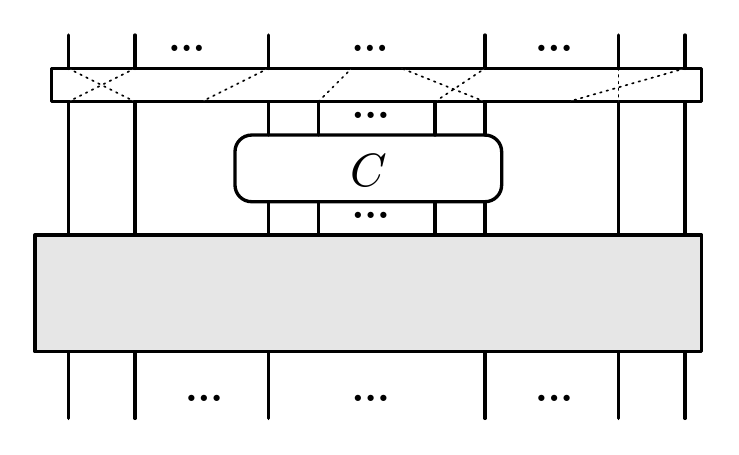}
  \caption{Circuits}\label{fig:normalFormCircuit}
\end{figure}
By Lemma~\ref{lem:bijectionStrictlyIncreasing}, we know that the cell~$C$ is not injective on inputs
of weight~$1$. It means there are two input wires~$i_1$ and~$i_2$ s.t.~$C$ gives the
same value on the two elements of~$\B^n$ consisting of~$0$s and a single~$1$
in position~$i_1$ or~$i_2$. Because this is independent of the inputs~$\vec v$
on the~$n-m$ remaining
wires, we obtain:
\begin{claim}
  Supposing contractive increasing boolean operators were finitely generated with
  a cells of arity less than~$m$, then for any non-invertible
  contractive increasing boolean operator~$f:\B^n\to\B^n$, with~$n\geq m$, we
  have:
  \[
    \exists \sigma\in\Perm_n \quad
    \forall \vec v \in \B^{n-m}
    \quad f\big(\sigma(01.\vec0.\vec v)\big) = f\big(\sigma(10.\vec0.\vec
    v)\big)
    \ .
  \]
\end{claim}
The permutation~$\sigma$ is used to simplify the notation: it reorders the
wires to put~$i_1$ and~$i_2$ in positions~$1$ and~$2$, and the remaining input
wires for~$C$ in positions~$3,\dots,c$.

For any maximal arity~$m$, we will construct a (large)~$n$ together with a
function~$f:\B^n\to\B^n$ that contradicts this fact: whenever we choose input
wires~$i_1$ and~$i_2$ and put any~$m-2$ other input wires to~$0$, we can
complete the remaining input wires in such a way that putting~$i_1\eqdef0$
and~$i_2\eqdef1$, or putting~$i_2\eqdef0$ and~$i_1\eqdef1$ makes a difference
in the output of the function. Thus, this function will not be representable
using the given set of cells.

\noindent
Given a (large)~$n$, define~$f:\B^n\to\B^n$ as:
\begingroup
\[
  f(u) \quad\eqdef\quad
  \cases{
    0^n & if $|u|=0$
            \hfill\sttt{(1)}\cr
    1\ 0^{n-1} & if $|u|=1$
            \hfill\sttt{(2)}\cr
    1^k\ 0^{n-k} & if $|u|=k$ is even
            \hfill\sttt{(3)}\cr
\noalign{\medbreak}
    1101\ 0^{n-4} & if $u = 0\cdots0\ 110^l1\ 0\cdots0$, with~$l>0$
            \hfill\sttt{(4)}\cr
    1110\ 0^{n-4} & if $|u| = 3$ but $u$ not of the previous shape
            \hfill\sttt{(5)}\cr
\noalign{\medbreak}
    1^{2^k}01\ 0^{n-2^k-2} &  if $u = 0\cdots0\ 1^{2^k}0^{2^k}1\ 0\cdots0$, with $k>1$
            \hfill\sttt{(6)}\cr
    1^{2^k}01\ 0^{n-2^k-2} &  if $u = 0\cdots0\ 10^{2^k}1^{2^k}\ 0\cdots0$, with $k>1$
            \hfill\sttt{(7)}\cr
    1^{2k}10\ 0^{n-2k-2} &  in all the remaining cases.
            \hfill\sttt{(8)}\cr
  }
\]
\endgroup
This function is contractive because we have~$|f(u)| = |u|$. Moreover, it is
increasing because whenever~$v$ is a successor\footnote{$v$ is a
\emph{successor} of~$u$ if~$v>u$ and~$|v|=|u|+1$.} of~$u$, we
have~$f(v)>f(u)$:
\begin{itemize}
  \item $f(u) = 1^{2k}\ 0\cdots$ when~$|u|=2k$
  \item $f(u) = 1^{2k}10\ 0\cdots$ or~$f(u) = 1^{2k}01\ 0\cdots$
    when~$|u|=2k+1$.
\end{itemize}
Suppose input wires~$k_1$, \dots,~$k_{m-2}$ are fixed to~$0$ and we want to
differentiate between input wires~$i_1$ and~$i_2$, with~$i_1<i_2$. By putting
some~$1$s in the appropriate remaining wires, we can make~$f$ give different results
when~``$i_1\eqdef0$, $i_2\eqdef1$'' and~``$i_1\eqdef1$, $i_2\eqdef0$''.

\begin{itemize}

  \item If there are two consecutive wires between~$i_1$ and~$i_2$ (but not
    touching~$i_2$) which are not among~$k_1$, \dots, $k_{m-2}$, we put those
    two wires to~$1$ and all the other wires to~$0$. By lines~\sttt{(4)}
    and~\sttt{(5)} from the definition of~$f$, we will get two different
    results.

  \item If not, the wires~$i_1$ and~$i_2$ cannot be too far apart. (There are at
    most~$2m-2$ wires between them...) If we can find a sequence of~$2^k$
    consecutive wires at distance~$2^k$ to the left of~$i_1$, or a sequence
    of~$2^k$ consecutive wires at distance~$2^k$ to the right of~$i_2$, we can
    put those wires to~$1$ and the rest to~$0$. By lines~\sttt{(6)}
    and~\sttt{(8)} or~\sttt{(7)} and~\sttt{(8)} of
    the definition of~$f$, we will also get different results.

\end{itemize}

For this to work, we have to make sure~$n$ is big enough. At worst, the
wires~$k_1$,\dots,~$k_{m-2}$ can prevent us from finding an appropriate
sequence~$m-2$ times. In particular, if~$i_1$ is big enough (bigger than
$2^{m+1}$), such a sequence is bound to happen. The same is true when~$i_2$ is
small enough compared to~$n$. In the end, choosing~$n$ bigger than,
say,~$2^{2m+2}$ plus an additional~$\varepsilon$ will guarantee that we can
differentiate any~$i_1$ and~$i_2$ among any set of~$m$ wires. A more careful
analysis shows that it is in fact enough to take~$n\eqdef2^{m+1}+4$. This
concludes the proof.
\qed
\end{pf}

%
%
%
%
%

\bibliographystyle{elsarticle-harv}


\bibliography{order}

\appendix

\section{An algorithm} 

Proposition~\ref{prop:characterization} is more elegant but
Lemma~\ref{lem:orderCharacterization} has an
interesting byproduct: it gives a concrete algorithm to check if~$\vec
X\sq\vec Y$.
For that, construct the function~$f$ from Lemma~\ref{lem:orderCharacterization} and check that it
satisfies the condition. Just as a proof of concept, here is the main part of
the algorithm, in the Python programming language. Minor alterations have been
made to make it more readable. The most difficult (and fun) part was to write
the function~\texttt{combinations} that generates all the vectors of
length~\texttt{n} and weight~\texttt{w} using one of the subtle algorithms
from~\cite{Knuth}!\footnote{The complete file is available from
\url{http://lama.univ-savoie.fr/~hyvernat/research.php}}

\lstset{
  language=Python,
  basicstyle=\scriptsize,
  tabsize=4,
  commentstyle=\ttfamily,
  }
\begin{lstlisting}
def check(N,n,X,Y):
	# N is a set, n is an integer, X / Y are tuples of sets.
	def combinations(w):
		# generates all vectors of weight w
        # omitted (see Knuth, or use you favorite method)
	def sup(u,v):           # complexity: O(n)
		# computes the pointwise "or" on n-tuples
		# omitted (simple)
	def weight(u):          # complexity: O(n)
		# computes the weight of an n-tuple
		# omitted (simple)
	def chi(a,Z):           # complexity: O(n log(z)) (z is cardinality of Z)
		for i in range(n):  # we use Python builtin "set" type
			if a in Z[i]:   # so that "a in Z[i]" mean "a belongs to Z[i]"
				u[i] = 1
		return u
	
	F = {}                  # F is a finite map with at most 2^n elements,
	                        # access is logarithmic: O(log(2^n)) = O(n)
	for a in N:                         # complexity: c *
	  chiX = chi(a,X)                   #                   n log(x)
	  chiY = chi(a,Y)                   #                 + n log(y)
	  F[chiY] = sup(F[chiY] , chiX)     #                 + 2n
	for w in range(n+1):                # generating all tuples
		for u in combinations(w):       # complexity : about 2^n *
			v = F[u]                                    #           n
			for i in range(n):                          #         + n *
				if u[i] == 1:                           #
					u[i] = 0                            #
					v = sup(v,F[u])                     #               n^2
					u[i] = 1                            #
			F[u] = v                                    #
			if weight(v) > w:                           #         + n
				return False
	return True     # if we reached this far, the condition is satisfied
\end{lstlisting}
If~$N$ has cardinality~$c$ and the components of~$\vec X$ and $\vec Y$ have
cardinalities at most~$x$ and~$y$; and if we suppose that the standard
operations on sets and finite functions have logarithmic complexity, the hints
in the comments give a total complexity of
roughly~$O\big(cn\big(\log(x)+\log(y)\big) + 2^n n^3\big)$.
Because both~$x=O(c)$ and~$y=O(c)$, we get a complexity of~$O\big(nc\log(c) +
2^n n^3\big)$. If~$c$ is fixed, this is~$O(n^3 2^n)$; if~$n$ is fixed, this
is~$O\big(c\log(c)\big)$.
In almost all cases, this is better (and much easier to write) than the naive
approach that checks if each~$\vec a\in\vec X$ is a section of~$\vec Y$, even
if we are allowed to use an oracle to guess the permutations.


\end{document}